\documentclass[12pt]{article}
\renewcommand{\baselinestretch}{1.2}

\newtheorem{thm}{Theorem}
\newtheorem{lem}[thm]{Lemma}

\newtheorem{prop}[thm]{Proposition}

\def\b/{$B_1^n$}
\def\bb/{$B_2^n$}
\def\elip/{ellipsoid}

\def\BM/{Brunn-Minkowski}
\def\bod/{convex body}
\def\bods/{convex bodies}
\def\d/{Dvoretzky}
\def\ep/{\varepsilon}
\def\iq/{inequality}
\def\iqs/{inequalities}
\def\iso/{isoperimetric}

\def\P/{probability}
\def\sph/{$S^{n-1}$}

\def\startproof{\bigskip \noindent {\em Proof }}
\def\endproof{\hfill $\vcenter{\hrule height .3mm
		\hbox {\vrule width .3mm height 2mm \kern 2mm
			\vrule width .3mm} \hrule height .3mm}$ \bigskip}

\title{Rational approximations to the zeta function II}
\author{Keith Ball}

\begin{document}

	\maketitle
	
	\begin{abstract}
		This note describes continued fraction representations for the rational approximations to $\zeta$ recently found by the author. It is tempting to think that these continued fractions might be analysed using a souped up version of the Worpitzky argument so as to produce zero-free regions for the approximations.
	\end{abstract}
	
	\section*{Introduction}
	The author's previous article \cite{rapp} describes 
	a sequence of rational functions which approximate $\zeta$ at least in the critical strip.
	
	The rational functions in question are the ratios
	\[ \frac{F_m(s)}{(s-1)G_m(s)} \]
	where for each $m$, $(s-1) G_m(s)$ is a rational function of $s$ that is close to $h_m^{1-s} \Gamma(s)$ and $F_m(s)$ is a rational function of $s$ that is close to $h_m^{1-s} \Gamma(s) \zeta(s)$ (and $h_m$ is the harmonic number $\sum_1^m 1/j$).

The sequence converges locally uniformly to $\zeta$, at least to the right of 
	the line $\{s: \Re s =0 \}$, (with the obvious convention at $s=1$). 
	After I circulated the article a number of people asked me whether my sequence could be generated by a simple continued fraction and I said I did not believe it. However  each individual rational function can of course be written as a continued fraction in many different ways. The aim of this article is to describe one way which produces a fraction that may be susceptible of analysis.
	
	By using the fact that $G_m$ approximates the Gamma function with an error that is quite easy to estimate one can check that $G_m$ does not vanish at $s=\sigma+i t$ if $0 < \sigma < 1$ and $|t|$ is smaller than a multiple of $\log m$. But the argument relies on the fact that $\Gamma$ itself has no zeros. Obviously if we are to use the functions $F_m$ to understand $\zeta$ we need an ``intrinsic'' way to find zero-free regions: an approach which uses only the very special shape of the rational functions.
	
	In Section \ref{Gamma} we shall see that the function $G_m$ has a continued fraction which can be directly analysed without reference to any property of the Gamma function so as to demonstrate that $G_m$ has no zeros $\sigma + i t$ in the critical strip where $|t|$ is smaller than about $\sqrt{\log m}$. In Section \ref{Zeta} we shall prove that the function $F_m$ has a continued fraction which shares at least some important features with the one for $G_m$.
	
	\section{An asymptotic series and a continued fraction for $G_m$}
	\label{Gamma}
	The rational approximations to $\zeta$ whose existence was proved in \cite{rapp} are defined as follows.
	For each integer $m \geq 0$ we define 
	\[ p_m(t)= \left(1-t \right) \left(1-\frac{t}{2} \right) \ldots \left( 1-\frac{t}{m} \right) \]
	and the coefficients $(a_{m,j})$ by
	\[ p_m(t)=\sum_0^m (-1)^j a_{m,j} t^j. \]
	We then set
	\[ F_m(s)=\sum_0^m \frac{a_{m,j} B_j}{s+j-1} \]
	where the $B_j$   are the usual Bernoulli numbers 
	and
	\[ G_m(s)=\sum_{j=0}^m (-1)^j \frac{a_{m,j}}{s+j-1}. \]

	For example
	\[ F_3(s)= \frac{1}{s-1} - \frac{11}{12 s} + \frac{1}{6 (s + 1)} = 
	\frac{3s^2+10s+11}{12 (s-1)s(s+1)} \]
	and 
		\[ G_3(s)= \frac{1}{s-1} - \frac{11}{6 s} + \frac{1}{s + 1}-\frac{1}{6(s+2)} = 
		\frac{s^2+6s+11}{3 (s-1)s(s+1)(s+2)}. \]

	The function $G_m$ is a rational function with poles at $1$, $0$, $-1$ and so on, which decays like $1/s^2$ at infinity. As a result we can expand it as a linear combination of reciprocals
	\[ \frac{1}{(s-1)s}, \; \; \; \frac{1}{(s-1)s(s+1)}, \; \; \; \frac{1}{(s-1)s(s+1)(s+2)} \ldots. \]
	The coefficients are essentially just the coefficients  $a_{m-1,j}$ defined above. More precisely
	\[ m s (s-1) G_m(s)= a_{m-1,0}+\frac{2 a_{m-1,1}}{s+1}+\frac{6 a_{m-1,2}}{(s+1)(s+2)}+\cdots. \]
		To see this observe that for $\Re s>1$
	\[ G_m(s)= \int_0^1 p_m(x) x^{s-2} \, dx. \]
	The behaviour of $G_m$ for large $s$ depends upon the behaviour of $p_m$ near $x=1$. But
	\begin{eqnarray*} G_m(s) & =& \int_0^1 p_m(1-x) (1-x)^{s-2} \, dx= \frac{1}{m} \int_0^1 x p_{m-1}(-x) (1-x)^{s-2} \, dx \\
	& = & \frac{1}{m} \int_0^1 \sum_0^{m-1} a_{m-1,j} \, x^{j+1}  (1-x)^{s-2} \, dx.
	\end{eqnarray*}
	
The sequence of coefficients $a_{m-1,j}$ is rather regular: they are all positive and the sequence is logarithmically concave because they are the coefficients of a polynomial $x \mapsto p_{m-1}(-x)$ whose zeros lie on the negative real axis. For small $j$ the coefficient $a_{m,j}$ is roughly $(\log m)^j/j!$. In particular $a_{m-1,1}$ is the harmonic number $h_{m-1}=1+1/2+1/3+\cdots+1/(m-1)$. 

On the face of it the formula

	\[ m s (s-1) G_m(s)= a_{m-1,0}+\frac{2 a_{m-1,1}}{s+1}+\frac{6 a_{m-1,2}}{(s+1)(s+2)}+\cdots \]
does not look as though it can tell us much about where $G_m$ is non-zero. Certainly if $s$ is large relative to $\log m$ then 
the sum is non-zero because it is dominated by the first term. But our interest is in finding zero-free regions that expand with $m$, rather than contracting: we want to understand what happens for $s$ {\em smaller} than say $\log m$.
However
a standard identity often known as Euler's continued fraction enables us to convert the ``asymptotic'' sum into a continued fraction for the reciprocal of $G_m$ (in which we drop the subcript $m-1$ from the coefficients for clarity)

\[\frac{1}{m s (s-1)G_m(s) }-1 \hspace{4in} \]
{\renewcommand{\baselinestretch}{1}
	\normalsize
	\[
	\begin{array}{rl}
	= \ {\strut \displaystyle -2 a_{1} \over\displaystyle 2 a_{1}+a_{0}(s+1)-
		{\strut \displaystyle 3 a_{0}a_{2}(s+1) \over\displaystyle 3 a_{2}+a_{1}(s+2)-
			{\strut \displaystyle 4 a_{1}a_{3}(s+2)  \over\displaystyle 4 a_{3}+a_{2}(s+3)}}}
	&  \\
	&  \ddots
	\end{array}  \] }
	
	There is a beautiful theorem of Worpitzky, see for example \cite{H} p.506, which shows that a continued fraction cannot ``blow up'' (cannot have zero denominator) if the denominators of the fraction are fairly large compared to the numerators. The hypothesis is that the product of two successive denominators should have absolute value at least 4 times as large as that of the numerator in between.

The point of this section is to observe that the fraction above representing $G_m$ has a structure which is well-adapted to Worpitzky's Theorem. The hypothesis in this case requires that for each $k$
\[ \left| ((k+1) a_{k}+a_{k-1}(s+k))((k+2) a_{k+1}+a_{k}(s+k+1)) \right| \geq 4 \left| (k+2) a_{k-1}a_{k+1}(s+k) \right|. \]
Set $v_k=\frac{(k+1) a_k}{(k+s) a_{k-1}}$ for each $k$.
Then the hypothesis is
\[ \left| (v_k+1)\left(1+\frac{1}{v_{k+1}} \right) \right| \geq 4. \]
The logarithmic concavity of the sequence $(a_k)$ shows that the sequence of ratios $a_k/a_{k-1}$ is decreasing. If $s$ is real and (say) in $(0,1)$ then $(v_k)$ is also decreasing and hence
\[ (v_k+1)\left(1+\frac{1}{v_{k+1}} \right) \geq (v_{k+1}+1)\left(1+\frac{1}{v_{k+1}} \right)=
  2+v_{k+1}+\frac{1}{v_{k+1}} \geq 4. \]
  
  Clearly in order to move off the real line one needs to understand how the strict logarithmic concavity of the sequence $(a_k/a_{k-1})$ is enough to compensate for the effect of the complex number $s$. The tricky point occurs where the sequence $a_k$ attains its maximum and $v_k$ is therefore close to 1. For larger values of $m$ this maximum occurs at values of $k$ close to $\log m$. At this point the ratio $\frac{(k+1)}{(k+s)}$ should have less effect when $m$ is large. So one hopes that the provable zero-free region should expand with $m$.
  
  In fact we can prove the following:
 \begin{prop} \label{prop} If $m$ is a natural number and for each $k$ we set $a_k=a_{m-1,k}$ and $v_k=\frac{(k+1) a_k}{(k+s) a_{k-1}}$ then we have
 \[ \left| (v_k+1)\left(1+\frac{1}{v_{k+1}} \right) \right| \geq 4 \]
 for $1 \leq k \leq m-2$ as long as $s=\sigma+i t$ satisfies $0<\sigma <1$ and $|t|< 1/2 \sqrt{\log m}$.
 
 \end{prop}

Newton's inequalities tell us that because $p_{m-1}(-x) = \sum_0^{m-1} a_j x^j$ has all its zeros on the negative real axis, the sequence 
\[ \frac{a_j}{{m-1 \choose j}} \]
is logarithmically concave.
This shows that the sequence 
\[ \left(  \frac{j}{m-j} \frac{a_j}{a_{j-1}} \right)_j \]
is decreasing and hence so is the sequence
\[ \left(  j \frac{a_j}{a_{j-1}} \right)_j. \]

Since $a_0=1$ and $a_1=h_{m-1}$ we can deduce that for each $j$
\[ \frac{a_j}{a_{j-1}} \leq \frac{h_{m-1}}{j}. \]
We will need an approximate reverse inequality which we prove as  a lemma.

\begin{lem} \label{triv} With the notation above 
	\[ \frac{a_j}{a_{j-1}} \geq \frac{h_{m-1}-1}{j} \]
	as long as $1 \leq j \leq h_{m-1}/2$.
	\end{lem}

\startproof The number $a_{j-1}$ is the sum of all products of $j-1$ distinct numbers in the set of reciprocals $1/r$ for $1 \leq r \leq m-1$. If we multiply this by the sum $h_{m-1} =\sum_1^{m-1} 1/r$ we obtain all possible products of $j$ distinct factors, each one repeated $j$ times, together with some products involving $j-2$ distinct factors and a squared factor.

Therefore
\[ h_{m-1} a_{j-1} \leq j a_j+a_{j-2} \sum_1^{m-1} \frac{1}{r^2} \leq j a_j+2 a_{j-2} . \]
If we write $w_{j}$ for the ratio $\frac{a_{j+1}}{a_j}$ we have
\[ \frac{h_{m-1}}{j} \leq w_{j-1}+\frac{2}{w_{j-2}}. \]
This inequality and a trivial induction show that
\[ w_{j-1} \geq \frac{h_{m-1}-1}{j} \]
as long as $1 \leq j \leq h_{m-1}/2$ as  required. \endproof

 \startproof (Of Proposition \ref{prop}) By the remarks above
\[ (k+1) \frac{a_{k+1}}{a_k} \leq k \frac{a_k}{a_{k-1}}. \]

Now $1+v_k= 1+ \frac{(k+1) a_k}{(k+s) a_{k-1}}$ and it easy to check that the absolute value of
\[ 1+\frac{(k+1)}{(k+s)} w \] increases with $w>0$ as long as $s$ lies in the critical strip.
Hence in proving the inequality we want, we may replace $\frac{a_k}{a_{k-1}}$ by  the smaller number 
$\frac{(k+1) a_{k+1}}{k a_{k}}$ or the still smaller number $\frac{(k+2) a_{k+1}}{(k+1) a_{k}}$.

So we then want to prove that 
\[ \left| \left( 1+\frac{(k+2) a_{k+1}}{(k+s) a_{k}} \right)
\left(1+\frac{(k+s+1) a_k}{(k+2) a_{k+1}} \right) \right| \geq 4. \]
If we write $w_k$ for the ratio $\frac{a_{k+1}}{a_k}$ (as before) this inequality becomes
\[ \left| 1+\frac{(k+s+1)}{(k+s)}+
\frac{k+2}{k+s} w_k +\frac{k+s+1}{k+2} \frac{1}{w_k} \right| \geq 4. \]
The expression inside the absolute value is
\[  2+\frac{1}{(k+s)}+
\frac{k+2}{k+s} w_k +\frac{k+s+1}{k+2} \frac{1}{w_k}. \]
Each of the terms in this sum has positive real part as long as $s$ is in the critical strip so it suffices to show that for each $k$
\begin{equation} \label{ineq} \Re \left(
\frac{k+2}{k+s} w_k +\frac{k+s+1}{k+2} \frac{1}{w_k} \right) \geq 2. \end{equation}

To handle the delicate range of $k$ observe that 
\[ \Re \frac{k+s+1}{k+2} = \frac{k+1+\sigma}{k+2} \]
and it easy to check that if $t^2 \leq k$ then
\[ \Re \frac{k+2}{k+s} \geq \frac{k+2}{k+1+\sigma}. \]
In this case, if we set $\theta$ to be the positive real number $\frac{k+2}{k+1+\sigma}$, we have 
\[ \Re \left(
\frac{k+2}{k+s} w_k +\frac{k+s+1}{k+2} \frac{1}{w_k} \right) \geq \theta w_k+\frac{1}{\theta w_k} \geq 2. \]

So it only remains to check (the trivial case) that (\ref{ineq}) is true for values of $k$ smaller than $1/4 \log m$ and 
$|t| < 1/2 \sqrt{\log m}$. By Lemma \ref{triv}
\[ w_k \geq \frac{h_m-1}{k+1} \]
and this implies that 
\[ \Re \frac{k+2}{k+s} w_k \geq 2. \]

\endproof
  
  \section{A continued fraction for $\zeta$}
  \label{Zeta}
  
  In the case of $F_m$ the approximation picks up the trivial zeros of $\zeta$ at $-2,-4,\ldots$ as far as $1-m$. Equivalently the function $F_m$ has poles at $1$, $0$, $-1$, $-3$ and so on but not at the even negative integers. As a result we can express $F_m$ as a sum
	\[ F_m(s)= \frac{b_0}{s-1}+\frac{b_1}{(s-1)s}+\frac{b_2}{(s-1)s(s+1)}+ \frac{b_3}{(s-1)s(s+1)(s+3)}+\cdots \]
	where after the first two terms we only use factors $s+2j-1$ to increase the degree of the denominator. It follows from remarks in \cite{rapp} that $b_0=1/(m+1)$ and hence that we can rewrite the sum as
	\[ (m+1) s F_m(s)= c_{m,0}+\frac{c_{m,1}}{s-1}+\frac{2 c_{m,2}}{(s-1)(s+1)}+\cdots\]
	where $c_{m,0}=1$.
	
From the definition of $F_m$,
\[ F_m(s)= \sum_0^m \frac{a_{m,j} B_j}{s+j-1}, \]
we get that the residue of $s F_m(s)$ at $1-j$ is $a_{m,j}(1-j) B_j$ and hence that
\[ (m+1) s F_m(s)-1= (m+1) \sum_0^m \frac{a_{m,j}(1-j) B_j}{s+j-1}. \]
The sum only involves even values of $j$ so we may write
\begin{eqnarray*} (m+1) s F_m(s)-1 & = & (m+1) \sum_{j \leq m/2} \frac{a_{m,2 j}(1-2 j) B_{2j}}{s+2 j-1} \\
& = &  \frac{(m+1)}{2} \int_0^1 \sum_{j \leq m/2} a_{m,2 j} (1-2j) B_{2 j} x^{j}  x^{(s-1)/2-1} \, dx \\
& = &  \frac{(m+1)}{2} \int_0^1 \sum_{j \leq m/2} a_{m,2 j} (1-2j) B_{2 j} (1-x)^{j}  (1-x)^{(s-1)/2-1} \, dx. \end{eqnarray*}
For $k>0$ let \begin{equation}
\label{ceq} c_{m,k}=(m+1)(-1)^{k-1} \sum_{j \leq m/2} a_{m,2 j} (1-2j) B_{2 j} {j \choose k-1}. \end{equation}
Then
\begin{eqnarray*} (m+1) s F_m(s)-1 & = & \frac{1}{2} \int_0^1 \sum_{k \leq m/2} c_{m,k+1} x^k  (1-x)^{(s-1)/2-1} \, dx \\
& = & \frac{c_{m,1}}{s-1}+\frac{2 c_{m,2}}{(s-1)(s+1)}+\cdots+ \frac{2^{j-1} c_{m,j} (j-1)!}{(s-1)\ldots(s+2j-3)}+\cdots. 
\end{eqnarray*}
	
Numerical evidence suggests that the coefficients $c_{m,j}$ for $F_m$ have similar properties to the $a_{m,j}$: for example $c_{m,1}=2(m+1)/(m+2) h_{m+1}$, the next coefficient $c_{m,2}$ grows like $(\log m)^2$ and so on. However it is not clear from the expression (\ref{ceq}) even that the coefficients are all positive. This will be demonstrated below.

This series for $F_m(s)$ can be converted into a continued fraction for $1/F_m$ much like the one for $1/G_m$:
\[\frac{1}{(m+1) s F_m(s) }-1 \hspace{4in} \]
{\renewcommand{\baselinestretch}{1}
	\normalsize
	\[
	\begin{array}{rl}
	= \ {\strut \displaystyle -c_{1} \over\displaystyle c_{1}+c_{0}(s-1)-
		{\strut \displaystyle 2 c_{0}c_{2}(s-1) \over\displaystyle 2 c_{2}+c_{1}(s+1)-
			{\strut \displaystyle 4 c_{1}c_{3}(s+1)  \over\displaystyle 4 c_{3}+c_{2}(s+3)}}}
	&  \\
	&  \ddots
	\end{array}  \] }
It is tempting to wonder whether the Worpitzky argument by itself gives non-trivial zero-free regions for $\zeta$ but my feeling is that it will not: that we will need a more subtle way to handle the expression for $1/F_m$ than we needed for $1/G_m$. It {\em does} appear to be the case that the ratio $c_k/c_{k-1}$ is decreasing. If this is true it would indicate that this representation for $F_m$ lies ``at the edge'' of what we need to prove zero-free regions. The stronger statement used above for the $a_k$ clearly cannot hold and numerically one can find values of $m$ for which it does not: for which the sequence $k c_k/c_{k-1}$ is {\em not} decreasing.

The expression for $c_{m,k}$ given in (\ref{ceq}) is not easy to understand directly: the alternation of sign in the Bernoulli numbers creates a subtle cancellation between the terms. However it is possible to prove that the coefficients are all positive.
To begin with we shall find a generating function.

\begin{lem} \label{genlem} For each $m$ and $k$ larger than 1 set
\[ c_{m,k}=(m+1)(-1)^{k-1} \sum_{j \leq m/2} a_{m,2 j} (1-2j) B_{2 j} {j \choose k-1}. \]
Then for $|z|, |y|<1$ and using principal values for the square root and logarithm,
\[ 1+\frac{1}{m+1} \sum_{k \geq 1,m\geq 1} c_{m,k} z^{k-1} y^m = (\log(1-y))^2 \frac{\partial}{\partial y} \frac{\sqrt{1-z}}{(1-y)^{\sqrt{1-z}}-1}. \]
\end{lem}

Once this lemma is established we can prove positivity using a standard continued fraction. The series for $(\log(1-y))^2$ has only non-negative coefficients so it suffices to check that the coefficients are positive in the expansion of  
\[ \frac{\partial}{\partial y} \frac{\sqrt{1-z}}{(1-y)^{\sqrt{1-z}}-1}. \]
It is known how to expand $(1-y)^x$ as a continued fraction (see for example \cite{H} p.535).
\[ \frac{(1-y)^t-1}{t}\hspace{4in} \]
{\renewcommand{\baselinestretch}{1}
	\normalsize
\[	\begin{array}{rl}
	  = \ {\strut \displaystyle 2 y \over\displaystyle 2-y+t y-
		{\strut \displaystyle (1-t^2)y^2 \over\displaystyle 3(2-y)-
			{\strut \displaystyle (4-t^2)y^2 \over\displaystyle 5(2-y)}}}
	&  \\
	&  \ddots
	\end{array}  \] }
	From this it follows that
	\[\frac{2 \sqrt{1-z}}{(1-y)^{\sqrt{1-z}}-1} \hspace{4in} \]
{\renewcommand{\baselinestretch}{1}
	\normalsize
	\[
	\begin{array}{rl}
	= -\frac{2}{y}+1-\sqrt{1-z}+\ {\strut \displaystyle z y \over\displaystyle 3(2-y)-
		{\strut \displaystyle (3+z)y^2 \over\displaystyle 5(2-y)-
			{\strut \displaystyle (8+z)y^2 \over\displaystyle 7(2-y)}}}
	&  \\
	&  \ddots
	\end{array}  \] }
	
When this expression is differentiated with respect to $y$ the first term gives $2/y^2$ which has a positive coefficient and the next term disappears. So it suffices to check that if the continued fraction is expanded as a power series in $y$ the coefficients are polynomials in $z$ with positive coefficients. If we want to check the coefficient of $y^m$ we only need to use the first $m/2$ levels of the continued fraction. Now start at the bottom of this finite continued fraction and work back up inductively. At each stage you have a fraction of the form
\[ \frac{(k^2-1+z)y^2}{(2k+1)(2-y)-b_2 y^2-b_3 y^3-\cdots} \]
where each $b_j$ is a polynomial in $z$ with positive coefficients. When this expression is expanded as  a power series in $y$
the coefficients are again polynomials in $z$ with positive coefficients.

Now for the proof of Lemma \ref{genlem}.

\startproof 
For $y$ and $z$ sufficiently small
\[ \frac{\sqrt{1-z}}{(1-y)^{\sqrt{1-z}}-1}= \sum_{j \geq 0} B_{j} (\sqrt{1-z})^{j} \frac{(\log(1-y))^{j-1}}{j!}, \]
its derivative  with respect to $y$ is
\[ \sum_{j \geq 0}  (1-j) B_{j} (\sqrt{1-z})^{j} \frac{1}{1-y}\frac{(\log(1-y))^{j-2}}{j!} \]
and the series is absolutely convergent.
Therefore
\[ (\log(1-y))^2 \frac{\partial}{\partial y} \frac{\sqrt{1-z}}{(1-y)^{\sqrt{1-z}}-1}= \sum_{j \geq 0}  (1-j) B_{j} (\sqrt{1-z})^{j} \frac{1}{1-y}\frac{(\log(1-y))^{j}}{j!}.\]
Since only the terms for which $j$ is even contribute to the sum we can introduce a negative sign to get
\begin{equation}
\label{expand} \sum_{j \geq 0}  (1-j) B_{j} (\sqrt{1-z})^{j} \frac{1}{1-y}\frac{(-\log(1-y))^{j}}{j!}. \end{equation}

It is a standard property of Stirling numbers that for each $j$ and for $|y|<1$
\[ \frac{(-\log(1-y))^{j+1}}{(j+1)!} = \sum_{m \geq 0} a_{m,j} \frac{y^{m+1}}{m+1} \]
and hence that
\[ \sum_{m \geq 0} a_{m,j} y^m= \frac{d}{dy} \frac{(-\log(1-y))^{j+1}}{(j+1)!}=\frac{1}{1-y}\frac{(-\log(1-y))^{j}}{j!}. \]

So the expression in (\ref{expand}) is equal to
\[ \sum_{j \geq 0} \sum_{m \geq 0} a_{m,j} y^m (1-j) B_{j} (\sqrt{1-z})^{j}. \]
Now replace $j$ by $2j$ using the fact that only even numbered terms occur to get 
\[ \sum_{j \geq 0} \sum_{m \geq 0} a_{m,2j} y^m (1 -2 j) B_{2 j} (1-z)^{j}. \]
If $m=0$ then $a_{m,j}$ is non-zero only if $j=0$ so the expression is
\[ 1+\sum_{j \geq 0} \sum_{m \geq 1} a_{m,2j} y^m (1 -2 j) B_{2 j} (1-z)^{j}. \]
Since we know that the series converges even if we replace the Bernoulli numbers by their absolute values we know that this series is absolutely convergent so we may interchange the order of summation to get
\[ 1+\sum_{m \geq 1} y^m \sum_{0 \leq j \leq m/2} a_{m,2j} (1 -2 j) B_{2 j} (1-z)^{j} \hspace{2in} \]
\begin{eqnarray*}  & = & 1+ \sum_{m \geq 1} y^m \sum_{0 \leq j  \leq m/2} \sum_{k=1}^{j+1} (-1)^{k-1}  a_{m,2 j} (1-2j) B_{2 j} {j \choose k-1} z^{k-1} \\ & = & 1+\frac{1}{m+1} \sum_{m \geq 1} y^m \sum_{k \geq 1} c_{m,k} z^{k-1}. 
\end{eqnarray*} 
\endproof

Once one has seen the Worpitzky argument and the generating function for the coefficients $c_{m,k}$ one is tempted to replace the functions $F_m$ with an analogous family of approximations indexed by the variable $y$ rather than by the power $m$ of $y$. To be precise we choose a large positive number $r$ and then for $\Re s >1$
\[ s r^{1-s}\zeta(s) \Gamma(s) = \frac{1}{4}\int_0^{\infty} \frac{r^2}{\sinh^2(r x/2)} x^s \, dx. \]
Now approximate the function by truncating the integral at $x=1$ and then substitute $x =\sqrt{1-z}$
to get
\[ \frac{1}{8}\int_0^{1} \frac{r^2 (1-z)}{\sinh^2(r \sqrt{1-z}/2)} (1-z)^{(s-1)/2-1} \, dz. \]
Now for each fixed $r$ we expand the function 
\[ \frac{r^2 (1-z)}{\sinh^2(r \sqrt{1-z}/2)} \]
as a power series in $z$ to obtain coefficients that replace the $c_{m,k}$.
Very limited numerical experiments suggest that this coefficient sequence has the ``right shape'' for each $r$.
The coefficient sequence is certainly logrithmically concave as one can check by using the Hadamard product for $\sinh$ and the Brunn-Minkowski inequality.

The resulting approximations to $\zeta$ don't have the appealingly simple matrix representations of the $F_m$ discussed in the first article in this series but they appear to have simpler coefficients in the ``asymptotic'' picture discussed in this article. My guess is that ultimately this simplicity is an illusion, but perhaps not. What is easy to see is that these coefficiants form a logarithmically concave sequence. The Hadamard product formula for $\sinh$ shows us that
\[ \frac{r^2 (1-z)}{\sinh^2(r \sqrt{1-z}/2)}=\frac{1}{(1+r^2(1-z)/(4 \pi^2))^2}\frac{1}{(1+r^2(1-z)/(16 \pi^2))^2}\ldots . \]
It is easily checked by hand that each factor in this product has a logarithmically concave coefficient sequence. By the discrete form of the Brunn-Minkowski inequality the product does as well.

	\end{document}